\documentclass{ita}
\newcommand{\D}{\mathbb{D}}
\newcommand{\PP}{\mathcal{P}}
\newcommand{\E}{\text{\rm Env}}
\newcommand{\w}{\omega}
\newcommand{\T}{\Theta}

\begin{document}
\title{Envelope Words and the Reflexivity of the Return Word Sequences in the Period-doubling Sequence}

\thanks{The research is supported by the Grant NSFC No.11431007 and No.11371210.}
%
\author{Yu-Ke HUANG}
\address{School of Mathematics and Systems Science, Beihang University (BUAA), Beijing, 100191, P. R. China; e-mail: huangyuke@buaa.edu.cn (Corresponding author).}
\author{Zhi-Ying WEN}
\address{Department of Mathematical Sciences, Tsinghua University, Beijing, 100084, P. R. China; e-mail: wenzy@tsinghua.edu.cn.}
%
\date{}
\begin{abstract}
We consider the infinite one-sided sequence over alphabet $\{a,b\}$ generated by the period-doubling substitution
$\sigma(a)=ab$ and $\sigma(b)=aa$, denoted by $\D$.
Let $r_p(\w)$ be the $p$-th return word of factor $\w$.
The main result of this paper is twofold. (1) For any factor $\w$ in $\D$, the return word sequence $\{r_p(\w)\}_{p\geq1}$ is $\T_1$ or $\T_2$. Both of them are substitutive sequences and determined completely in this paper.
(2) For any factor $\w$ in $\T_1$ (resp. $\T_2$), the return word sequence $\{r_p(\w)\}_{p\geq1}$ is still $\T_1$ or $\T_2$.
We call it the reflexivity property of the return word sequence.
As an application, we introduce a notion of spectrum for studying some typical combinatorial properties, such as separated, adjacent and overlapped.
\end{abstract}

\maketitle
\section*{Introduction}

Let $\mathcal{A}=\{a,b\}$ be a binary alphabet.
The period-doubling sequence $\D$ is the fixed point beginning with the letter $a$ of substitution $\sigma(a)=ab$ and $\sigma(b)=aa$,
and also the first difference of the Thue-Morse sequence. Here we use an equivalent substitution $\hat{\sigma}(1)=10$, $\hat{\sigma}(0)=11$,
and the definition of the difference of an integer sequence is natural.
It has been heavily studied in mathematics and computer science.
Damanik \cite{D2000} determined the numbers of palindromes, squares and cubes of length $n$ occurring in $\D$.
Allouche-Peyri\`ere-Wen-Wen~\cite{APWW1998}
proved that all the Hankel determinants of $\D$ are odd integers.
See also \cite{ACF2006,B2006,CLR2015,CRS2012,D1998-2,D2001,FH2015,GW2014,L1998,PRRV2015,RSV2013}.

\medskip

We denote $A_m=\sigma^m(a)$ and $B_m=\sigma^m(b)$ for $m\geq0$.
Then $|A_m|=|B_m|=2^m$, where $|\omega|$ is the length of $\w$.
Let $\delta_m\in\{a,b\}$ be the last letter of $A_m$. Obviously,
$\delta_m=a$ if and only if $m$ even; $\delta_{m+1}$ is the last letter of $B_m$.
For later use, we list some elementary properties of $\D$, which can be proved easily by induction:
$A_{m}\delta_m^{-1}=B_{m}\delta_{m+1}^{-1}$ and
$A_{m}=a\prod_{j=0}^{m-1} B_{j}$ for $m\geq0$.
Moreover $\D=a\prod_{j=0}^\infty B_j$.

The mirror word of $\omega=x_1x_2\cdots x_n$ is defined to be $\overleftarrow{\omega}=x_n\cdots x_2x_1$. A word $\omega$ is called a palindrome if $\omega=\overleftarrow{\omega}$.
We define $\overline{\cdot}:\mathcal{A}\rightarrow\mathcal{A}$ by $\overline{a}=b$ and $\overline{b}=a$.
We say that $\nu$ is a prefix (resp. suffix) of a word $\w$ if there exists word $u$ such that $\w=\nu u$ (resp. $\w=u\nu$), $|u|\geq0$, which denoted by $\nu\triangleleft\w$ (resp. $\nu\triangleright\w$).
In this case, we write $\nu^{-1}\w=u$ (resp. $\w\nu^{-1}=u$), where $\nu^{-1}$ is the inverse word of $\nu$ such that $\nu\nu^{-1}=\nu^{-1}\nu=\varepsilon$ (empty word).

Let $\rho=x_1x_2\cdots x_n$ be a finite word.
For any $1\leq i\leq j\leq n$,
we define
$\rho[i]=\rho[i,i]=x_i$, $\rho[i,i-1]=\varepsilon$,
$\rho[i,j]=x_ix_{i+1}\cdots x_{j-1}x_j$ and
$\rho[i,\mathrm{end}]=\rho[i,n]$.
Since $\D$ is uniformly recurrent, each factor $\w$ appears infinite many times in the sequence, which is arranged as $\w_p$ $(p\ge 1)$. Here we distinguish $\w_p\neq\w_q$ if $p\neq q$.
Let $[\w]_p$ be the position of the first letter of $\w_p$. We also call it the position of $\w_p$ for short.
The definition of return word is from Durand \cite{D1998}.
Let $\w\prec\D$ be a factor of $\D$.
The \emph{$p$-th return word} of $\w$ is $r_{p}(\w)=\D[[\w]_p,[\w]_{p+1}-1]$ for $p\geq1$.
The sequence $\{r_p(\w)\}_{p\geq1}$ is called the
\emph{return word sequence} of factor $\w$.
We denote $r_0(\w)=\D[1,[\w]_1-1]$ the prefix of $\D$ before $\w_1$, which is not a return word.

\medskip

Huang-Wen \cite{HW2015-1,HW2015-2} determine the structure of the return word sequences of the Fibonacci sequence $\mathbb{F}$ and the Tribonacci sequence $\mathbb{T}$, respectively.
More precisely, for any factor $\omega$ of $\mathbb{F}$ (resp. $\mathbb{T}$),
the return word sequence $\{r_p(\omega)\}_{p\geq1}$ is still $\mathbb{F}$ (resp. $\mathbb{T}$).
Using these properties, we determined the numbers of palindromes, squares and cubes occurring in $\mathbb{F}[1,n]$ (the prefix of $\mathbb{F}$ of length $n$) and $\mathbb{T}[1,n]$ for all $n\geq1$.
These topics are of great importance in computer science.

The main tool of the two papers above is ``kernel word''.
However, in the study of the period-doubling sequence $\mathbb{D}$, the technique of kernel word fail. In fact, there is no kernel set in $\mathbb{D}$, which satisfies the ``uniqueness of kernel decomposition'' property, see \cite{HW2015-1,HW2015-2}.
To overcome this difficulty, we introduce a new notion called ``envelope word'' and prove the ``uniqueness of envelope extension'' property.

The other difficulty is: the return word sequence of $\D$ is not unique, see examples in Section 1. In this paper, we determine the return word sequences completely. We also find out the relations among them and $\D$, called the reflexivity property.

\section{Main Results}

Now we turn to list the main results, and give some examples.

Let $\{E^i_{m}\mid i=1,2,\ m\geq1\}$ be a set of factors with
$$E^1_{m}=A_{m}\delta_{m}^{-1}\text{ and }E^2_{m}=B_{m}B_{m-1}\delta_{m}^{-1}.$$
We call $E^i_{m}$ the $m$-th \emph{envelope word} of type $i$.

Obviously, the lengths of all envelope words are odd, except $E^2_1=aa$.
\begin{center}
\begin{tabular}{|l|l|l|l|l|l|}
\multicolumn{6}{c}{Tab.1: The first few values of $A_m$, $B_m$, $E^1_{m}$ and $E^2_{m}$}\\\hline
$m$&0&1&2&3&4\\\hline
$A_m$&$a$&$ab$&$abaa$&$abaaabab$&$abaaabababaaabaa$\\\hline
$B_m$&$b$&$aa$&$abab$&$abaaabaa$&$abaaabababaaabab$\\\hline
$E^1_{m}$&&$a$&$aba$&$abaaaba$&$abaaabababaaaba$\\\hline
$E^2_{m}$&&$aa$&$ababa$&$abaaabaaaba$&$abaaabababaaabababaaaba$\\\hline
\end{tabular}
\end{center}

Define $\T_1=\tau_1(\D)$ and $\T_2=\tau_2(\D)$, where $\tau_1(a)=a$, $\tau_1(b)=bb$, $\tau_2(a)=ab$ and $\tau_2(b)=acac$.
Obviously, $\T_1$ and $\T_2$ are $\D$ over the alphabets $\{a,bb\}$ and
$\{ab,acac\}$, respectively.
\begin{center}
\begin{tabular}{|lc|l|}
\multicolumn{3}{c}{Tab.2: The first few letters of $\D$, $\T_1$ and $\T_2$}\\\hline
$\D$&=&  $abaaabababaaabaaabaaabababaaabababaaabababaaabaaabaaab\cdots$\\\hline
$\T_1$&=&$abbaaabbabbabbaaabbaaabbaaabbabbabbaaabbabbabbaaabbabbab\cdots$\\\hline
$\T_2$&=&$abacacabababacacabacacabacacabababacacabababacacabababac\cdots$\\\hline
\end{tabular}
\end{center}

\textbf{First result.} For any factor $\w$ of $\D$, denoted by $\w\prec\D$,
the return word sequence $\{r_p(\w)\}_{p\geq1}$ is $\Theta_1$ or $\Theta_2$.
For instance, $\{r_p(aba)\}_{p\geq1}$ in $\D$ is $\Theta_1$ over alphabet $\{r_1(aba),r_2(aba)\}=\{abaa,ab\}$, which is denoted by $\{A,B\}$ as below.
\begin{equation*}
\begin{split}
\D=&\underbrace{abaa}_A\underbrace{ab}_B\underbrace{ab}_B
\underbrace{abaa}_A\underbrace{abaa}_A\underbrace{abaa}_A\underbrace{ab}_B\underbrace{ab}_B
\underbrace{abaa}_A\underbrace{ab}_B\underbrace{ab}_B\\
&\underbrace{abaa}_A\underbrace{ab}_B\underbrace{ab}_B
\underbrace{abaa}_A\underbrace{abaa}_A\underbrace{abaa}_A\underbrace{ab}_B\underbrace{ab}_B
\underbrace{abaa}_A\underbrace{abaa}_A\underbrace{abaa}_A\ \cdots
\end{split}
\end{equation*}
And $\{r_p(aa)\}_{p\geq1}$ in $\D$ is $\Theta_2$ over
$\{r_1(aa),r_2(aa),r_4(aa)\}=\{a,aababab,aab\}$, which is denoted by $\{\alpha,\beta,\gamma\}$ as below.
\begin{equation*}
\begin{split}
\D=&ab\underbrace{a}_\alpha\underbrace{aababab}_\beta \underbrace{a}_\alpha \underbrace{aab}_\gamma \underbrace{a}_\alpha \underbrace{aab}_\gamma \underbrace{a}_\alpha \underbrace{aababab}_\beta \underbrace{a}_\alpha\\
&\underbrace{aababab}_\beta\underbrace{a}_\alpha\underbrace{aababab}_\beta \underbrace{a}_\alpha
\underbrace{aab}_\gamma \underbrace{a}_\alpha\underbrace{aab}_\gamma \underbrace{a}_\alpha\underbrace{aababab}_\beta\ \cdots
\end{split}
\end{equation*}

\textbf{Second result.} For any factor $\w$ of $\T_1$ (resp. $\T_2$), the return word sequence $\{r_p(\w)\}_{p\geq1}$ is $\T_1$ or $\T_2$ too.
For instance:

1. take $\w=abba\prec \T_1$, \hspace{0.08cm}then $\{r_p(\w)\}_{p\geq1}$ in $\T_1$ is $\T_1$ over $\{r_1(\w),r_2(\w)\}$;

2. take $\w=aa\prec \T_1$,  \hspace{0.34cm} then $\{r_p(\w)\}_{p\geq1}$ in $\T_1$ is $\T_2$ over $\{r_1(\w),r_2(\w),r_4(\w)\}$;

3. take $\w=aba\prec \T_2$,  \hspace{0.14cm} then $\{r_p(\w)\}_{p\geq1}$ in $\T_2$ is $\T_1$ over $\{r_1(\w),r_2(\w)\}$;

4. take $\w=ababa\prec \T_2$, then $\{r_p(\w)\}_{p\geq1}$ in $\T_2$ is $\T_2$ over $\{r_1(\w),r_2(\w),r_4(\w)\}$.


\smallskip

Generally speaking, Let $\tau$ be a sequence and $f$ be a mapping on $\tau$, where $f(\tau)$ maps $\tau$ to the set $\big\{\{r_p(\omega)\}_{p\geq1}\mid \omega\prec\tau\big\}$.
An immediate corollary of this result is the set $\{\Theta_1,\Theta_2\}$ is invariant under the mapping $f$, i.e., $f(\Theta_i)=\{\Theta_1,\Theta_2\}$ for $i=1,2$.
We call it the reflexivity property of the return word sequences.
Recall that for any factor, the return word sequence of the Fibonacci (resp. Tribonacci) sequence is itself \cite{HW2015-1,HW2015-2}. They are special cases of the reflexivity  property.

\section{The return word sequences of envelope words in $\D$}

Our goal in this section is to prove two theorems below. They determine the return word sequences for all envelope words.

\begin{thrm}\label{T2.1}\
The return word sequence $\{r_p(E^1_{m})\}_{p\geq1}$ in $\D$ is $\Theta_1$ over the binary alphabet
$\{A_{m},A_{m-1}\}$.
\end{thrm}

\begin{thrm}\label{T2.2}\
The return word sequence $\{r_p(E^2_{m})\}_{p\geq1}$ in $\D$  is $\Theta_2$ over the three-letter alphabet
$\{A_{m-1},A_{m-1}A_{m}B_{m+1},B_{m}B_{m-1}\}$.
\end{thrm}

\subsection{Proof of Theorem \ref{T2.1}}

We first give a criterion to determine all occurrences of a factor, called \emph{Criterion~E}.
It is useful in our proofs. Let $\w=u\nu$ where $|u|,|\nu|>0$. Obviously, if $\D[L,L+|\omega|-1]=\omega$, then $\D[L,L+|u|-1]=u$.
Thus in order to determine all occurrences of $\w$, we need two steps:

1. find out all occurrences of $u$;

2. check every occurrences of $u$. When an occurrence of $u$ is followed by $\nu$, we say the factor $u$ at this position can extend to $\w$.

\begin{lmm}\label{L2.1}\
$A_m$ occurs exactly twice in $A_mA_m$ (resp. $A_mB_mA_m$) for $m\geq0$.
\end{lmm}

\begin{proof} The result is clearly true for $m=0$.
Assume the result is true for $m$. Then all possible positions of
$A_{m}$ in $A_{m+1}A_{m+1}$ and $A_{m+1}B_{m+1}A_{m+1}$
are shown with ``underbrace'' as below.
\begin{equation*}
\begin{split}
A_{m+1}A_{m+1}&=\underbrace{A_{m}}_{[1]}B_{m}\underbrace{A_{m}}_{[2]}B_{m},\\
A_{m+1}B_{m+1}A_{m+1}&=\underbrace{A_{m}}_{[3]}B_{m}\underbrace{A_{m}}_{[4]} \underbrace{A_{m}}_{[5]}\underbrace{A_{m}}_{[6]}B_{m}.
\end{split}
\end{equation*}
Among them, only the $A_{m}$ at positions [1], [2], [3] and [6] are followed by $B_m$.
By Criterion E, $A_{m+1}=A_mB_m$ occurs twice in $A_{m+1}A_{m+1}$ (resp. $A_{m+1}B_{m+1}A_{m+1}$).
So the conclusion hold for $m\geq0$ by induction.
\end{proof}

By Lemma \ref{L2.1}, the first three positions of $A_{m-1}$ are shown below.
$$\underbrace{A_{m-1}}_{[1]}B_{m-1}\underbrace{A_{m-1}}_{[2]}
\underbrace{A_{m-1}}_{[3]}A_{m-1}\triangleleft A_{m+2}\triangleleft\D.$$
All of them are followed by $B_{m-1}\delta_{m}^{-1}=A_{m-1}\delta_{m-1}^{-1}$.
By Criterion E, all of them can extend to $E^1_{m}=A_{m-1}B_{m-1}\delta_{m}^{-1}$.
Thus the first three positions of $E^1_{m}$ are:
$$[E^1_{m}]_1=1,~[E^1_{m}]_2=2^m+1,~[E^1_{m}]_3=3\times2^{m-1}+1.$$
By the definition of return words $r_p(E^1_{m})$, we have
\begin{equation*}
\begin{cases}
r_0(E^1_{m})=\D[1,[E^1_{m}]_1-1]=\D[1,0]=\varepsilon;\\
r_1(E^1_{m})=\D[[E^1_{m}]_1,[E^1_{m}]_2-1]=\D[1,2^m]=A_{m};\\
r_2(E^1_{m})=\D[[E^1_{m}]_2,[E^1_{m}]_3-1]=\D[2^m+1,3\times2^{m-1}]=A_{m-1}.
\end{cases}
\end{equation*}
Moreover $|r_0(E^1_{m})|=0$, $|r_1(E^1_{m})|=2^{m}$ and $|r_2(E^1_{m})|=2^{m-1}$.

\medskip

For $m\geq1$ we define an infinite sequence of alphabets
$$\mathcal{A}^1_{m}=\{\mathbf{a},\mathbf{b}\}=\{r_1(E^1_{m}),r_2(E^1_{m})r_2(E^1_{m})\}
=\{A_{m},B_{m}\}.$$
For each fixed $m\geq1$, we denote factors
$\sigma^n(\mathbf{a})$ and $\sigma^n(\mathbf{b})$ over alphabet $\mathcal{A}^1_{m}$ by
$\mathbf{A_n}$ and $\mathbf{B_n}$, respectively.
For instance, taking $m=1$, $\mathcal{A}^1_{1}=\{A_{1},B_{1}\}=\{ab,aa\}$.
Then $\mathbf{A_2}=\sigma^2(\mathbf{a})=\sigma^2(ab)=abaaabab=A_3$ over alphabet $\{a,b\}$.

\begin{lmm}\label{L2.2}
Over the alphabet $\mathcal{A}^1_{m}$, 
$\mathbf{A_n}=A_{m+n}$, $\mathbf{B_n}=B_{m+n}$ for $n\geq0$.
\end{lmm}

\begin{proof} 
Since $\mathbf{A_0}=\mathbf{a}=A_{m}$, $\mathbf{B_0}=\mathbf{b}=B_{m}$, the result is correct for $n=0$. Assume it is true for $n$, we prove it for $n+1$.
\begin{equation*}
\begin{split}
\mathbf{A_{n+1}}&=\mathbf{A_nB_n}=A_{m+n}B_{m+n}=A_{m+n+1},\\
\mathbf{B_{n+1}}&=\mathbf{A_nA_n}=A_{m+n}A_{m+n}=B_{m+n+1}.
\end{split}
\end{equation*}
By induction, the result is correct for all $n\geq0$.
\end{proof}

By the definitions of $\mathcal{A}^1_{m}$ and $\Theta_1$, Theorem \ref{T2.1} is equivalent to
Theorem \ref{T2.1}'.

\medskip

\noindent\textbf{Theorem \ref{T2.1}'.}
\emph{For each envelope word $E^1_{m}$ ($m\geq1$), the return word sequence $\{r_p(E^1_{m})\}_{p\geq1}$ in the period-doubling sequence is itself still a period-doubling sequence over the alphabet $\mathcal{A}^1_{m}$.}

\begin{proof} By $A_{m}=a\prod_{j=0}^{m-1}B_j$ and Lemma \ref{L2.2}, we have
\begin{equation*}
\begin{split}
\D=a\prod_{j=0}^\infty B_j=a\prod_{j=0}^{m-1}B_j\prod_{j=m}^\infty B_j
=\mathbf{a}\prod_{j=0}^\infty B_{m+j}=\mathbf{a}\prod_{j=0}^\infty \mathbf{B_{j}}
=\mathbf{D}.
\end{split}
\end{equation*}
Here $\D$ and $\mathbf{D}$ are the period-doubling sequence over alphabets
$\mathcal{A}$ and $\mathcal{A}^1_{m}$, respectively.
Thus the conclusion holds.
\end{proof}

\subsection{Proof of Theorem \ref{T2.2}}

By an analogous argument, for $m\geq1$, $B_m$ occurs

(1) once in $A_mB_m$, at position $2^{m}+1$;

(2) three times in $B_mA_mB_m$, at positions 1, $2^{m-1}+1$ and $2^{m+1}+1$;

(3) three times in $B_mA_mA_mA_mB_m$, at positions 1, $2^{m-1}+1$ and $2^{m+2}+1$.

Thus by $A_mB_mA_mA_mA_mB_mA_mB_mA_mB_m\triangleleft A_{m+4}\triangleleft\D$ and Criterion E, the first five positions of $E^2_{m}=B_mB_{m-1}\delta_m^{-1}$ are:
$$2^m+1,~3\times 2^{m-1}+1,~5\times 2^m+1,~11\times2^{m-1}+1,~7\times 2^m+1.$$
The expressions of $r_p(E^2_{m})$, $m\geq1$ and $1\leq p\leq4$ are
\begin{equation*}
\begin{cases}
r_0(E^2_{m})=A_{m},&|r_0(E^2_{m})|=2^{m};\\
r_1(E^2_{m})=r_3(E^2_{m})=A_{m-1},&|r_1(E^2_{m})|=2^{m-1};\\
r_2(E^2_{m})=A_{m-1}A_{m}B_{m+1},&|r_2(E^2_{m})|=7\times2^{m-1};\\
r_4(E^2_{m})=B_{m}B_{m-1},&|r_4(E^2_{m})|=3\times2^{m-1}.
\end{cases}
\end{equation*}

For $m\geq1$ we define an infinite sequence of alphabets
$$\mathcal{A}^2_{m}=\{\mathbf{a},\mathbf{b}\}
=\{r_1r_2,r_1r_4r_1r_4\}
=\{A^{-1}_{m}A_{m+2}A_{m},A^{-1}_{m}B_{m+2}A_{m}\}.$$
Here we denote $r_i(E^2_{m})$ by $r_i$ for short ($i=1,2,4$).

For each fixed $m\geq1$, we denote factors
$\sigma^n(\mathbf{a})$ and $\sigma^n(\mathbf{b})$ over alphabet $\mathcal{A}^2_{m}$ by $\mathbf{A_n}$
and $\mathbf{B_n}$, respectively.
Similarly to Lemma \ref{L2.2}, we have Lemma \ref{L2.3}.
By the definitions of $\mathcal{A}^2_{m}$ and $\Theta_2$, Theorem \ref{T2.2} is equivalent to
Theorem \ref{T2.2}'.

\begin{lmm}\label{L2.3}
Over $\mathcal{A}^2_{m}$,
$\mathbf{A_n}=A^{-1}_{m}A_{m+n+2}A_{m}$, $\mathbf{B_n}=A^{-1}_{m}B_{m+n+2}A_{m}$
for $n\geq0$.
\end{lmm}

\noindent\textbf{Theorem \ref{T2.2}'.}
\emph{For each envelope word $E^2_{m}$ ($m\geq1$), the return word sequence $\{r_p(E^2_{m})\}_{p\geq1}$ in the period-doubling sequence is itself still a period-doubling sequence over the alphabet $\mathcal{A}^2_{m}$.}

\begin{proof} Denote $r_0(E^2_{m})$ by $r_0$. By $A_{m}=a\prod_{j=0}^{m-1} B_{j}$ and Lemma \ref{L2.3},
\begin{equation*}
\begin{split}
&\D=a\prod_{j=0}^\infty B_j=a\prod_{j=0}^{m+1}B_j\prod_{j=m+2}^\infty B_j
=A_{m+2}\prod_{j=0}^\infty B_{m+j+2}\\
=&A_{m} A^{-1}_{m}A_{m+2}A_{m} \prod_{j=0}^\infty (A^{-1}_{m}B_{m+j+2}A_{m})
=r_0\mathbf{a}\prod_{j=0}^\infty \mathbf{B_{j}}
=r_0\mathbf{D}.
\end{split}
\end{equation*}
Here $\D$ and $\mathbf{D}$ are the period-doubling sequence over alphabets
$\mathcal{A}$ and $\mathcal{A}^2_{m}$, respectively.
Notice that, $r_0$ is not a return word, so
we omit $r_0$ in the return word sequence. Thus the conclusion holds.
\end{proof}

\section{Uniqueness of envelope extension}

In this section, we give two types of \emph{uniqueness of envelope extensions}, which play important roles in our studies.
Using them, we can extend Theorems \ref{T2.1}, \ref{T2.2} and other related properties from envelope words to general factors.

First we define the \emph{order} of envelope words that
$E^1_{m}\sqsubset E^2_{m}$ and
$E^i_{m}\sqsubset E^j_{m+1}$ for $i,j\in\{1,2\}$, $m\geq1$.
For any factor $\w$, let
$$\E(\w)=\min_{\sqsubset}\{E^i_{m}\mid\w\prec E^i_{m},~i=1,2,~m\geq1\},$$
which is called \emph{the envelope of factor} $\w$.

By the definition of envelope, the envelope of each factor is unique.
But we do not know:
(1) whether $\w$ occurs in $\E(\w)$ only once or not;
(2) the relation between the positions of the $p$-th occurrence of factor $\w$ and the $p$-th occurrence of its envelope $\E(\w)$ for all $p\geq1$.

The two types of uniqueness of envelope extensions will answer these questions.

\begin{thrm}[The weak type of envelope extension]\label{T3.1}\
Each factor $\w$ occurs exactly once in $\E(\w)$.
\end{thrm}

Thus for each factor $\w$, there exist a unique integer $i$ ($0\leq i\leq |\E(\w)|-|\w|$), such that
\begin{equation}\label{E1}
\w=\E(\w)[i+1,i+|\w|].
\end{equation}

\begin{thrm}[The strong type of envelope extension]\label{T3.2}
Every occurrence of factor $\w$ in $\D$ can extend to a factor $\E(\w)$.
\end{thrm}

Thus for each factor $\w$ with expression~(\ref{E1}) and $p\geq1$, the difference $[\w]_p-[\E(\w)]_p$ is equal to $i$, depending only on $\w$.
This means, the $p$-th occurrence of factor $\w$ and the $p$-th occurrence of its envelope $\E(\w)$ are ``relative rest".
Theorem \ref{T3.2} is clearly true for $m=1,2$. For $m\geq3$, 
we prove it in two cases, see Propositions \ref{P3.9} and \ref{P3.11}.

\subsection{Basic properties of envelope words}

For later use, we give some basic properties of envelope words first.

Since $E^1_{m}=A_m\delta_{m}^{-1}$, $E^2_{m}=B_{m}B_{m-1}\delta_{m}^{-1}$ and $A_m\delta_{m}^{-1}=B_m\delta_{m+1}^{-1}$, we have
$$E^1_{m+1}=E^1_{m}\delta_{m}E^1_{m}\text{ and } E^2_{m+1}=E^1_{m}\delta_{m}E^1_{m}\delta_{m}E^1_{m}\text{ for }m\geq1.$$
Thus both $E^1_{m}$ and $E^2_{m}$ are palindromes.
By induction, we give a more general form in Lemma \ref{L3.1}.

\begin{lmm}[]\label{L3.1}
For $m>n$ and $i\in\{1,2\}$,
$E^i_{m}=E^1_{n}x_1E^1_{n}\cdots E^1_{n}x_hE^1_{n}$,
where $x_1\cdots x_h=\overline{E^i_{m-n}}$ for $n$ odd, and
$x_1\cdots x_h=E^i_{m-n}$ for $n$ even.
\end{lmm}

For instance, when $n=1$ (odd) and $m=3$, we give all $E^1_{1}=a$ with underline: $E^1_{3}=\underline{a}b\underline{a}a\underline{a}b\underline{a}$.
In this case, $x_1\cdots x_h=bab=\overline{aba}=\overline{E^1_{2}}=\overline{E^1_{m-n}}$.

\begin{lmm}[]\label{L3.2}
The factor $E^1_{n}\delta_{m}E^1_{k}$ ($n>k$) is a palindrome if and only if
$n-k=1$, $m$ and $k$ have the same parity.
In this case, $E^1_{n}\delta_{m}E^1_{k}=E^1_{n}\delta_{n-1}E^1_{n-1}=E^2_{n}$.
\end{lmm}

\begin{proof} (1) When $n=k+1$,
$E^1_{n}\delta_{m}E^1_{k}=E^1_{k}\underline{\delta_{k}}E^1_{k}\underline{\delta_{m}}E^1_{k}$.
Comparing the two letters with underlines, $E^1_{n}\delta_{m}E^1_{k}$
is a palindrome if and only if $\delta_{k}=\delta_{m}$. This means $m$ and $k$ have the same parity.

(2) When $n\geq k+2$,
$E^1_{n}\delta_{m}E^1_{k}=E^1_{k}x_1E^1_{k}x_2\cdots x_hE^1_{k}\delta_{m}E^1_{k}$.
By Lemma~\ref{L3.1}, $x_1x_2\cdots x_h=E^1_{n-k}$ or $\overline{E^1_{n-k}}$. Since
$n-k\geq 2$, $E^1_{n-k}[2]=b$ and $E^1_{n-k}[\mathrm{end}]=a$.
Thus $x_2\neq x_h$, $E^1_{n}\delta_{m}E^1_{k}$ can not be a palindrome.
So the conclusion holds.
\end{proof}

\begin{prpstn}[]\label{P3.1}
Let palindrome $\omega$ be a proper prefix of $E^1_{m}$ (resp. $E^2_{m}$). Then there exists $n\leq m$ such that $\omega=E^1_{n}$.
\end{prpstn}

\begin{proof} Case 1. Consider $E^1_{m}$ first.
We prove the result by induction on $m$. Clearly it is correct for $m=1,2$.
Now assume the result is true for $m$, we prove it for $m+1$. Since $E^1_{m+1}=E^1_{m}\delta_{m}E^1_{m}$, there are three cases.

(1) If $\w$ is a prefix of $E^1_{m}$, there exists $n\leq m$ such that $\omega=E^1_{n}$.

(2) For $m\geq2$, $E^1_{m}[2]=b\neq E^1_{m}[\mathrm{end}]=a$. So
$\w=E^1_{m}\delta_{m}=abE^1_{m}[3,\mathrm{end}-1]a\delta_{m}$ can not be a palindrome.

(3) If $\w=E^1_{m}\delta_{m}\mu$ is a palindrome, where $\mu$ is a proper prefix of $E^1_{m}$.
Since $\overleftarrow{\w}=\w$, we have $\w=\overleftarrow{\mu}E^1_{m}[|\mu|+1,\mathrm{end}]\delta_{m}\mu$.
This means $\overleftarrow{\mu}$ is a prefix of $E^1_{m}$ too.
So $\mu$ is palindrome too.
Thus there exists $n'< m$ such that $\w=E^1_{n'}$. This means $\w=E^1_{m}\delta_{m}E^1_{n'}$ is a palindrome. This contradicts Lemma \ref{L3.2}.

By induction, the conclusion holds for $E^1_{m}$, $m\geq1$.

\smallskip

Case 2. The conclusion for $E^2_{m}$ can be verified by induction too. We give only the induction step. Since $E^2_{m+1}=E^1_{m}\delta_{m}E^1_{m}\delta_{m}E^1_{m}$, there are seven cases as below, where $\mu$ is a proper prefix of $E^1_{m}$.
\begin{equation*}
\begin{cases}
(1)\w=\mu,~(2)\w=E^1_{m},~(3)\w=E^1_{m}\delta_m,\\
(4)\w=E^1_{m}\delta_m\mu,
~(5)\w=E^1_{m}\delta_mE^1_{m}=E^1_{m+1},\\
(6)\w=E^1_{m}\delta_mE^1_{m}\delta_m=E^1_{m+1}\delta_m
,~(7)\w=E^1_{m}\delta_mE^1_{m}\delta_m\mu=E^1_{m+1}\delta_m\mu.
\end{cases}
\end{equation*}
Obviously $\w$ is palindrome in cases (2) and (5). By the proof above, there exists $n\leq m$ such that $\omega=E^1_{n}$ in case (1); and $\w$ can not be palindrome in cases (3), (4), (6) and (7).
This completes the proof.
\end{proof}

\subsection{The simplification by palindromic property}

The ``weak type of envelope extension'' means ``each $\w$ occurs in $\E(\w)$ only once''.
Though the analysis in this subsection, we only need to prove the last property for $\w$ is a palindrome.
Moreover, if $\w$ occurs in $\E(\w)$ at least twice, we can pick two of them.
So we only need to negate the proposition that ``there exist a palindrome $\w$ occurs in $\E(\w)$ twice''.

\begin{lmm}[]\label{L3.4}
For each palindrome $\Lambda$ of odd length, if there exists a factor $\w$ satisfying:
(1) $\w$ occurs in $\Lambda$ twice;
(2) both of the two occurrences of $\w$ contain the middle letter of $\Lambda$.
Then there exists a palindrome occurring in $\Lambda$ twice with symmetric positions.
\end{lmm}

\begin{proof} Denote these two occurrences of $\w$ by $\w'$ and $\w''$. Moreover
$$\w'=\Lambda[d_1+1,d_1+|\w|]\text{ and }\w''=\Lambda[d_2+1,d_2+|\w|],$$
where $0\leq d_1<d_2\leq|\Lambda|-|\w|$.

Case 1. $d_1=|\Lambda|-d_2-|\w|$. $\w'$ and $\w''$ occur in $\Lambda$ at symmetric positions.
Since $\Lambda$ is palindrome,
$\overleftarrow{\w''}=\overleftarrow{\Lambda}[|\Lambda|-d_2-|\w|+1,|\Lambda|-d_2]
=\Lambda[d_1+1,d_1+|\w|]=\w'$.
Since $\w'$ and $\w''$ are two occurrences of $\w$, we have $\overleftarrow{\w}=\w$.
This means the factor $\w$ is a palindrome occurring at symmetric positions.

Case 2. $d_1\neq|\Lambda|-d_2-|\w|$. Since $\Lambda$ is palindrome, we assume $d_1<|\Lambda|-d_2-|\w|$ without loss of generality.
By the expressions of $\w'$ and $\w''$,
$$\overleftarrow{\w'}=\Lambda[|\Lambda|-d_1-|\w|+1,|\Lambda|-d_1]\text{ and }\overleftarrow{\w''}=\Lambda[|\Lambda|-d_2-|\w|+1,|\Lambda|-d_2].$$
Since both $\w'$ and $\w''$ contain the middle letter of $\Lambda$,
$\w'$ and $\overleftarrow{\w''}$ are overlapped. Define a factor
$\mu'=\Lambda[d_1+1,|\Lambda|-d_2]$. Then $\w'\triangleleft\mu'$, $\overleftarrow{\w''}\triangleright\mu'$ and $|\mu'|<|\w'|+|\w''|=2|\w|$.
This means $\mu'$ is palindrome.
Similarly, $\mu''=\Lambda[d_2+1,|\Lambda|-d_1]$ is palindrome.

Moreover, $\mu'$ and $\mu''$ occur at symmetric positions.
So the conclusion holds.
\end{proof}

\subsection{The weak type of envelope extension}

\begin{lmm}\label{L3.5}
Let $\w$ be a factor of $\D$. (1) If there exists integer $m\geq3$ such that $\E(\w)=E^1_{m}$, then $E^1_{m-2}\prec\w$; Moreover, one element in set
$$\{\delta_{m}E^1_{m-2}\delta_{m-1},\delta_{m-1}E^1_{m-2}\delta_{m-1},\delta_{m-1}E^1_{m-2}\delta_{m}\}$$ is a factor of $\w$.
(2) If there exists integer $m\geq3$ such that $\E(\w)=E^2_{m}$, then $\delta_{m-1}E^1_{m-1}\delta_{m-1}\prec\w$.
(3) Factor $\w$ contains the middle letter of $\E(\w)$.
\end{lmm}

\begin{proof} Case 1. $\E(\w)=E^1_{m}$.
\begin{figure}[!ht]
\centering
\footnotesize
\setlength{\unitlength}{0.54mm}
\begin{picture}(160,28)
\put(0,14){$E^1_{m}=E^1_{m-2}\delta_{m}E^1_{m-3}\delta_{m-1}E^1_{m-3}\delta_{m-1}E^1_{m-3}\delta_{m-1}E^1_{m-3}\delta_{m}E^1_{m-2}$}
\scriptsize
\put(17,0){\line(0,1){18}}
\put(81,0){\line(0,1){5}}
\put(81,10){\line(0,1){8}}
\put(40,0){$E^1_{m-1}$}
\put(55,1){\vector(1,0){26}}
\put(39,1){\vector(-1,0){22}}
\put(39,6){\line(0,1){12}}
\put(109,6){\line(0,1){12}}
\put(60,6){$E^2_{m-2}$}
\put(76,7){\vector(1,0){33}}
\put(59,7){\vector(-1,0){20}}
\put(100,21){$E^2_{m-2}$}
\put(66,12){\line(0,1){12}}
\put(136,12){\line(0,1){12}}
\put(116,22){\vector(1,0){20}}
\put(99,22){\vector(-1,0){33}}
\put(117,28){$E^1_{m-1}$}
\put(94,12){\line(0,1){8}}
\put(94,27){\line(0,1){5}}
\put(158,12){\line(0,1){20}}
\put(132,30){\vector(1,0){26}}
\put(116,30){\vector(-1,0){22}}
\end{picture}
\normalsize
\caption{The envelope words $E^1_{m-1}$ and $E^2_{m-2}$ occur in $E^1_{m}$.\label{Fig:2}}
\end{figure}
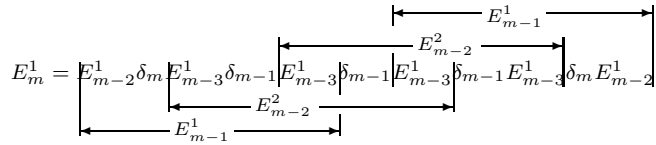

Since $\w$ can not be the factor of $E^1_{m-1}$ or $E^2_{m-2}$, $\w$ has three cases by Figure \ref{Fig:2}:

(i) \hspace{0.22cm} $\w\succ \delta_{m}E^1_{m-3}\delta_{m-1}E^1_{m-3}\underline{\delta_{m-1}}
=\delta_{m}E^1_{m-2}\delta_{m-1}$;

(ii) \hspace{0.11cm} $\w\succ \delta_{m-1}E^1_{m-3}\underline{\delta_{m-1}}E^1_{m-3}\delta_{m-1}
=\delta_{m-1}E^1_{m-2}\delta_{m-1}$;

(iii) \hspace{0.01cm} $\w\succ \underline{\delta_{m-1}}E^1_{m-3}\delta_{m-1}E^1_{m-3}\delta_{m}
=\delta_{m-1}E^1_{m-2}\delta_{m}$.

In each case, (1) $\w\succ E^1_{m-2}$; (2) $\w$ contains the middle letter $\delta_{m-1}$ of $\E(\w)$, which is shown with underline.

Case 2. $\E(\w)=E^2_{m}$.
\begin{figure}[!ht]
\centering
\setlength{\unitlength}{0.6mm}
\begin{picture}(90,18)
\put(0,7){$E^2_{m}=E^1_{m-1}\delta_{m-1}E^1_{m-1}\delta_{m-1}E^1_{m-1}$}
\small
\put(34,0){$E^1_{m}$}
\put(17,0){\line(0,1){11}}
\put(59,0){\line(0,1){11}}
\put(32,1){\vector(-1,0){15}}
\put(44,1){\vector(1,0){15}}
\put(62,14){$E^1_{m}$}
\put(44.5,6){\line(0,1){11}}
\put(87,6){\line(0,1){11}}
\put(60,15){\vector(-1,0){15.5}}
\put(73,15){\vector(1,0){14}}
\end{picture}
\normalsize
\caption{The envelope word $E^1_{m}$ occurs in $E^2_{m}$.\label{Fig:3}}
\end{figure}
Since $\w$ can not be the factor of $E^1_{m}$, by Figure \ref{Fig:3} we have
$\w\succ \delta_{m-1}E^1_{m-1}\delta_{m-1}$.
Obviously, $\w$ contains the middle letter of $\E(\w)$
\end{proof}

We call factor $u$ occurs in the middle of $\w$ if there exists integer $i$ ($0\leq i\leq|\w|$)
such that $u=\w[i+1,|\w|-i]$.


\begin{proof}[Proof of Theorem \ref{T3.1}]
The result is clearly true for $m=1,2$.

For $m\geq3$,
suppose the factor $\w$ occurs twice in $\E(\w)$, denoted by $\w'$ and $\w''$. Denote $\E(\w)$ by $E^i_{m}$ for $i\in\{1,2\}$.
Since $|E^i_{m}|$ is odd, by Lemmas \ref{L3.4} and \ref{L3.5}(3), we can assume without loss of generality that:
$\w$ is palindrome; $\w'$ and $\w''$ occur in $E^i_{m}$ at symmetric positions.

\smallskip

Case 1. $\E(\w)=E^1_{m}$.

$E^1_{m-2}$ occurs five times in $E^1_{m}$, denoted by positions [1] to [5] in Figure \ref{Fig:4}. By Lemma \ref{L3.5}(1), $\w'$ and $\w''$ contain at least one $E^1_{m-2}$ at positions [2] to [4].
\begin{figure}[!ht]
\centering
\footnotesize
\setlength{\unitlength}{0.50mm}
\begin{picture}(170,23)
\put(0,14){$E^1_{m}=\underbrace{E^1_{m-2}}_{[1]}
\delta_{m}E^1_{m-3}\delta_{m-1}E^1_{m-3}\delta_{m-1}E^1_{m-3}\delta_{m-1}E^1_{m-3}\delta_{m}
\underbrace{E^1_{m-2}}_{[5]}$}
\put(42,12){$\underbrace{\hspace{2.3cm}}_{[2]}$}
\put(103,12){$\underbrace{\hspace{2.3cm}}_{[4]}$}
\put(72,18){$\overbrace{\hspace{2.3cm}}$}
\tiny
\put(92,24){$[3]$}
\end{picture}
\normalsize
\caption{The five occurrences of $E^1_{m-2}$ in $E^1_{m}$.\label{Fig:4}}
\end{figure}

1. If $\w'$ and $\w''$ contain only one occurrence of $E^1_{m-2}$. Since they are palindrome, $E^1_{m-2}$ occurs in the middle of $\w'$ (resp. $\w''$). Since $\w'$ and $\w''$ occur in $E^1_{m}$ at symmetric positions, 
they contain $E^1_{m-2}$ at positions [2] and [4] respectively.
Since $E^1_{m-2}$ at [2] is preceded by $\delta_m$ and followed by $\delta_{m-1}$, $\w'$ can not be a palindrome.

2. If $\w'$ and $\w''$ contain two occurrences of $E^1_{m-2}$. Since they are palindrome, the two occurrences of $E^1_{m-2}$ occurs at symmetric positions.
There are two subcases: (i) $\w'$ contains $E^1_{m-2}$ at [1] and [2]. Then $\w'=E^1_{m-2}\delta_{m}E^1_{m-2}=E^1_{m-1}$, contradicting to $\E(\w)=E^1_{m}$. (ii) $\w'$ contains $E^1_{m-2}$ at [2] and [3]. Since $E^1_{m-2}$ is preceded by $\delta_m$ at [2], and followed by $\delta_{m-1}$ at [3], $\w'$ can not be a palindrome.

3. Similarly, when $\w'$ and $\w''$ contain three (resp. four, five) occurrences of $E^1_{m-2}$, we find contradictions.

\smallskip

Case 2. $\E(\w)=E^2_{m}$.

$E^1_{m-1}$ occurs three times in $E^2_{m}$, denoted by positions [1] to [3] as below. By Lemma \ref{L3.5}, $\w'$ and $\w''$ contain the $\delta_{m-1}E^1_{m-1}\delta_{m-1}$ in the middle of $E^2_{m}$.
$$E^2_{m}=\underbrace{E^1_{m-1}}_{[1]}\delta_{m-1}\underbrace{E^1_{m-1}}_{[2]}\delta_{m-1}\underbrace{E^1_{m-1}}_{[3]}$$
Obviously there are no palindromes $\w'$ and $\w''$ satisfy the assumptions.

Thus, each factor $\w$ occurs exactly once in $\E(\w)$.
\end{proof}

\subsection{The strong type of envelope extension}

In this subsection, we prove that
every occurrence of $\w$ in $\D$ can extend to $\E(\w)$. This means let factor $\w$ have expression~(\ref{E1}), then every occurrence of $\w$ in $\D$ must preceded by $\E(\w)[1,i]$ and followed by $\E(\w)[i+|\w|+1,\mathrm{end}]$.

\medskip

\textbf{Case 1.} $\E(\w)=E^1_{m}$ for $m\geq3$.

\medskip

By Theorem \ref{T2.1}, the return word sequence
$\{r_p(E^1_{m-2})\}_{p\geq1}$ is $\Theta_1$ over alphabet 
$\{A_{m-2},A_{m-3}\}
=\{E^1_{m-2}\delta_{m},E^1_{m-2}\delta_{m}B_{m-3}^{-1}\},$
where $E^1_{m-2}=A_{m-2}\delta_{m}^{-1}$.
Since $\Theta_1=\tau_1(\D)$ where $\tau_1(a)=a$ and $\tau_1(b)=bb$,
all factors in $\Theta_1$ of length 3 are $\{abb,bba,baa,aaa,aab,bab\}$.
Thus all occurrences of $E^1_{m-2}$ in $\D$ can be divided into 6 cases as below, where
$\mathbf{a}=r_1(E^1_{m-2})$, $\mathbf{b}=r_2(E^1_{m-2})$.
\begin{equation*}
\begin{cases}
\mathbf{abb}= A_{m-2} E^1_{m-2}\delta_{m}B_{m-3}^{-1} A_{m-3}\succ\delta_{m}E^1_{m-2}\delta_{m-1};\\
\mathbf{bba}= A_{m-3} E^1_{m-2}\delta_{m}B_{m-3}^{-1} A_{m-2}\succ\delta_{m-1}E^1_{m-2}\delta_{m-1};\\
\mathbf{baa}= A_{m-3} E^1_{m-2}\delta_{m} A_{m-2}\succ\delta_{m-1}E^1_{m-2}\delta_{m};\\
\mathbf{aaa}= A_{m-2} E^1_{m-2}\delta_{m} A_{m-2}\succ\delta_{m}E^1_{m-2}\delta_{m};\\
\mathbf{aab}= A_{m-2} E^1_{m-2}\delta_{m} A_{m-3}\succ\delta_{m}E^1_{m-2}\delta_{m};\\
\mathbf{bab}= A_{m-3} E^1_{m-2}\delta_{m} A_{m-3}\succ\delta_{m-1}E^1_{m-2}\delta_{m}.
\end{cases}
\end{equation*}
Here we use $\delta_{m}B_{m-3}^{-1} A_{m-3}=(B_{m-3}\delta_{m}^{-1})^{-1}(A_{m-3}\delta_{m-1}^{-1})\delta_{m-1}=\delta_{m-1}$,
and always rewrite the middle letter $\mathbf{a}$ or $\mathbf{b}$ to be the expression with envelope word $E^1_{m-2}$.

On the other hand, by Lemma \ref{L3.5}(1), one element in set
$$\{\delta_{m}E^1_{m-2}\delta_{m-1},\delta_{m-1}E^1_{m-2}\delta_{m-1},\delta_{m-1}E^1_{m-2}\delta_{m}\}$$
is a factor of $\w$. So we consider three subcases in Lemmas \ref{L3.6} to \ref{L3.8}.

\begin{lmm}\label{L3.6}
The factor $\delta_{m}E^1_{m-2}\delta_{m-1}$ occurs in the period-doubling sequence always preceded by $A_{m-2}\delta_{m}^{-1}$
and followed by $A_{m-1}\delta_{m-1}^{-1}$.
This means every occurrence of factor $\delta_{m}E^1_{m-2}\delta_{m-1}$ can extend to $E^1_{m}$.
\end{lmm}

\begin{proof} The factor $\delta_{m}E^1_{m-2}\delta_{m-1}$ occurs in the period-doubling sequence has only one case: $\mathbf{abb}$. Since $\Theta_1=\tau_1(\D)$ where $\tau_1(a)=a$ and $\tau_1(b)=bb$, each occurrence of $\mathbf{abb}$ can extend to $\mathbf{\underline{abb}abb}$ or $\mathbf{\underline{abb}aa}$ in $\Theta_1$.
\begin{equation*}
\begin{cases}
\mathbf{\underline{abb}abb}&=A_{m-2}E^1_{m-2}\delta_{m-1}A_{m-2}A_{m-3}A_{m-3}=A_{m-2}E^1_{m-2}\delta_{m-1}A_{m-1};\\
\mathbf{\underline{abb}aa}&=A_{m-2}E^1_{m-2}\delta_{m-1}A_{m-2}A_{m-2}
=A_{m-2}E^1_{m-2}\delta_{m-1}B_{m-1}.
\end{cases}
\end{equation*}

Since $B_{m-1}\delta_{m}^{-1}=A_{m-1}\delta_{m-1}^{-1}$,
both of them have prefix
$$A_{m-2}\delta_{m}^{-1}\underline{\delta_{m}E^1_{m-2}\delta_{m-1}} A_{m-1}\delta_{m-1}^{-1}
=A_{m}\delta_{m}^{-1}=E^1_{m}.$$
Here we show the factor $\delta_{m}E^1_{m-2}\delta_{m-1}$ with underline.
So the conclusion holds.
\end{proof}

\begin{lmm}\label{L3.7}
The factor $\delta_{m-1}E^1_{m-2}\delta_{m-1}$ occurs in the period-doubling sequence always preceded by the word $A_{m-2}A_{m-3}\delta_{m-1}^{-1}$
and followed by $B_{m-3}A_{m-2}\delta_{m}^{-1}$.
This means every occurrence of factor $\delta_{m-1}E^1_{m-2}\delta_{m-1}$ can extend to $E^1_{m}$.
\end{lmm}

\begin{proof} The factor $\delta_{m-1}E^1_{m-2}\delta_{m-1}$ occurs in $\D$ has only one case: $\mathbf{bba}$. It extends to $\mathbf{a\underline{bba}a}$ or $\mathbf{a\underline{bba}bb}$ in $\Theta_1$.
\begin{equation*}
\begin{cases}
\mathbf{a\underline{bba}a}
&=A_{m-2}A_{m-3}E^1_{m-2}\delta_{m-1}B_{m-3}A_{m-2};\\
\mathbf{a\underline{bba}bb}
&=A_{m-2}A_{m-3}E^1_{m-2}\delta_{m-1}B_{m-3}A_{m-3}A_{m-3}.
\end{cases}
\end{equation*}
Both of them have prefix
$A_{m-2}A_{m-3}\delta_{m-1}^{-1}\underline{\delta_{m-1}E^1_{m-2}\delta_{m-1}} B_{m-3}A_{m-2}\delta_{m}^{-1}
=E^1_{m}$.

This completes the proof.
\end{proof}

\begin{lmm}\label{L3.8}
The factor $\delta_{m-1}E^1_{m-2}\delta_{m}$ occurs in the period-doubling sequence always preceded by $A_{m-1}\delta_{m-1}^{-1}$
and followed by $A_{m-2}\delta_{m}^{-1}$.
This means every occurrence of factor $\delta_{m-1}E^1_{m-2}\delta_{m}$ can extend to $E^1_{m}$.
\end{lmm}

\begin{proof} The factor $\delta_{m-1}E^1_{m-2}\delta_{m}$ occurs in $\D$ has two cases: $\mathbf{baa}$ and $\mathbf{bab}$.
They extend to $\mathbf{ab\underline{baa}}$ and $\mathbf{ab\underline{bab}b}$ in $\Theta_1$, respectively.
\begin{equation*}
\begin{cases}
\mathbf{ab\underline{baa}}&=A_{m-2}A_{m-3}A_{m-3}E^1_{m-2}\delta_{m}A_{m-2}
=A_{m-1}E^1_{m-2}\delta_{m}A_{m-2};\\
\mathbf{ab\underline{bab}b}&=A_{m-2}A_{m-3}A_{m-3}E^1_{m-2}\delta_{m}A_{m-3}A_{m-3}
=A_{m-1}E^1_{m-2}\delta_{m}B_{m-2}.
\end{cases}
\end{equation*}
Both of them have prefix
$A_{m-1}\delta_{m-1}^{-1}\underline{\delta_{m-1}E^1_{m-2}\delta_{m}} A_{m-2}\delta_{m}^{-1}
=E^1_{m}$.

Thus the conclusion holds.
\end{proof}

\begin{prpstn}[]\label{P3.9}
If there exists integer $m\geq3$ such that $\E(\w)=E^1_{m}$, then every occurrence of $\w$ in the period-doubling sequence can extend to $\E(\w)$.
\end{prpstn}

\textbf{Case 2.} $\E(\w)=E^2_{m}$ for $m\geq3$.

\medskip

By Theorem \ref{T2.1},
the return word sequence $\{r_p(E^1_{m-1})\}_{p\geq1}$ in $\D$ is $\Theta_1$ over alphabet 
$\{A_{m-1},A_{m-2}\}
=\{E^1_{m-1}\delta_{m-1},E^1_{m-1}\delta_{m-1}B_{m-2}^{-1}\},$
where $E^1_{m-1}=A_{m-1}\delta_{m-1}^{-1}$.
Since all factors with length 3 in $\Theta_1$ are $\{abb,bba,baa,aaa,aab,bab\}$,
all occurrences of $E^1_{m-1}$ in $\D$ can be divided into 6 cases as below, where
$\mathbf{a}=r_1(E^1_{m-1})$, $\mathbf{b}=r_2(E^1_{m-1})$.
\begin{equation*}
\begin{cases}
\mathbf{abb}= A_{m-1} E^1_{m-1}\delta_{m-1}B_{m-2}^{-1} A_{m-2}\succ\delta_{m-1}E^1_{m-1}\delta_{m};\\
\mathbf{bba}= A_{m-2} E^1_{m-1}\delta_{m-1}B_{m-2}^{-1} A_{m-1}\succ\delta_{m}E^1_{m-1}\delta_{m};\\
\mathbf{baa}= A_{m-2} E^1_{m-1}\delta_{m-1} A_{m-1}\succ\delta_{m}E^1_{m-1}\delta_{m-1};\\
\mathbf{aaa}= A_{m-1} E^1_{m-1}\delta_{m-1} A_{m-1}\succ\delta_{m-1}E^1_{m-1}\delta_{m-1};\\
\mathbf{aab}= A_{m-1} E^1_{m-1}\delta_{m-1} A_{m-2}\succ\delta_{m-1}E^1_{m-1}\delta_{m-1};\\
\mathbf{bab}= A_{m-2} E^1_{m-1}\delta_{m-1} A_{m-2}\succ\delta_{m}E^1_{m-1}\delta_{m-1}.
\end{cases}
\end{equation*}
Here we always rewrite the middle letter $\mathbf{a}$ or $\mathbf{b}$ to be an expression with envelope word $E^1_{m-1}$.

On the other hand, by Lemma \ref{L3.5}(2), $\delta_{m-1}E^1_{m-1}\delta_{m-1}\prec\w$.

\begin{lmm}\label{L3.10}
The factor $\delta_{m-1}E^1_{m-1}\delta_{m-1}$ occurs in the period-doubling sequence always preceded and followed by the word $A_{m-1}\delta_{m-1}^{-1}$.
This means every occurrence of factor $\delta_{m-1}E^1_{m-1}\delta_{m-1}$ can extend to $E^2_{m}$.
\end{lmm}

\begin{proof} By the analysis above, $\delta_{m-1}E^1_{m-1}\delta_{m-1}$ occurs in $\D$ has two cases: $\mathbf{aaa}$ and $\mathbf{aab}$.
Here
$\mathbf{\underline{aaa}}=A_{m-1}E^1_{m-1}\delta_{m-1}A_{m-1}$.
Moreover $\mathbf{aab}$ extends to $\mathbf{\underline{aab}b}$ in $\Theta_1$.
$$\mathbf{\underline{aab}b}=A_{m-1}E^1_{m-1}\delta_{m-1}A_{m-2}A_{m-2}=A_{m-1}E^1_{m-1}\delta_{m-1}B_{m-1}.$$
Both of them have prefix
$A_{m-1}\delta_{m-1}^{-1}\underline{\delta_{m-1}E^1_{m-1}\delta_{m-1}} A_{m-1}\delta_{m-1}^{-1}
=E^2_{m}$.

This means the conclusion holds.
\end{proof}

\begin{prpstn}[]\label{P3.11}
If there exists integer $m\geq3$ such that $\E(\w)=E^2_{m}$, then every occurrence of $\w$ in the period-doubling sequence can extend to $\E(\w)$.
\end{prpstn}

\section{The return word sequences of general factors in $\mathbb{D}$}

Let factor $\w$ have expression~(\ref{E1}). By the strong type of envelope extension (see Theorem \ref{T3.2}),
we give the expressions of $r_p(\w)$ that:
\begin{equation}\label{E2}
r_0(\w)=r_0(\mathrm{E})\mathrm{E}[1,i]
\text{ and }
r_p(\w)=r_p(\mathrm{E})[i+1,\mathrm{end}]\mathrm{E}[1,i]\text{ for }p\geq1.
\end{equation}
Here we denote $\E(\w)$ by $\mathrm{E}$ for short.
An immediate corollary is $|r_p(\w)|=|r_p(\E(\w))|$ for $p\geq1$, which will be used in Section 6.

Using expression~(\ref{E2}), we extend Theorems \ref{T2.1} and \ref{T2.2} to general factors as Theorems \ref{T4.1} as below. This is the first main result in our paper.

\begin{thrm}[]\label{T4.1}\
Let $\w$ be a factor. (1) If $\E(\w)=E^1_{m}$,
the return word sequence $\{r_p(\w)\}_{p\geq1}$ in $\D$ is $\Theta_1$ over the alphabet $\{r_1(\w),r_2(\w)\}$.
(2) If $\E(\w)=E^2_{m}$,
the return word sequence $\{r_p(\w)\}_{p\geq1}$ in $\D$ is $\Theta_2$ over the alphabet $\{r_1(\w),r_2(\w),r_4(\w)\}$.
\end{thrm}

\section{The reflexivity of the return word sequences}

In Huang-Wen \cite{HW2015-1,HW2015-2}, we showed that: for any factor, the return word sequence of the Fibonacci (resp. Tribonacci) sequence is itself.
In this section, we show that:
for any factor $\w$ in $\T_1$ (resp. $\T_2$), the return word sequence $\{r_p(\w)\}_{p\geq1}$ is still $\T_1$ or $\T_2$.
We call it the reflexivity property of the return word sequence.
This is the second main result in our paper.

Recall that $\Theta_i=\tau_i(\D)$,
where $i=1,2$, $\tau_1(a)=a$, $\tau_1(b)=bb$, $\tau_2(a)=ab$ and $\tau_2(b)=acac$.
The methods are similar with the proofs of $\D$.
Two types of ``uniqueness of envelope extension" still hold.
Thus let factor $\w$ have expression~(\ref{E1}), the return words $r_p(\w)$ and $r_p(\E(\w))$ have relation (\ref{E2}) too.

Thus we only give the definitions of envelope words and the return word sequences of each types of envelop words in $\T_1$ and $\T_2$ in this section.
We also give some examples here. More examples, see Section 1.

The envelope words of $\Theta_1$ are factors in set
$\{{}^1\!E^{j}_{m}\mid j=1,2,~m\geq1\}$ where
\begin{equation*}
\begin{cases}
{}^1\!E^{1}_{m}=\tau_1(E^1_{m})\text{ for }m\geq1;\\
{}^1\!E^{2}_{1}=b\text{ and }{}^1\!E^{2}_{m}=\tau_1(E^2_{m-1})\text{ for }m\geq2.
\end{cases}
\end{equation*}
The envelope words of $\Theta_2$ are factors in set
$\{{}^2\!E^{j}_{m}\mid j=1,2,~m\geq1\}$ where
\begin{equation*}
\begin{cases}
{}^2\!E^{1}_{1}=a\text{ and }{}^2\!E^{1}_{m}=\tau_2(E^1_{m-1})a\text{ for }m\geq2;\\
{}^2\!E^{2}_{1}=aca\text{ and }{}^2\!E^{2}_{m}=\tau_2(E^2_{m-1})a\text{ for }m\geq2.
\end{cases}
\end{equation*}

\begin{center}
\begin{tabular}{|l|l|l|l|l|}
\multicolumn{5}{c}{Tab.3: The first few values of $E^j_{m}$ and ${}^i\!E^{j}_{m}$ for $i,j=1,2$}\\\hline
$m$&1&2&3&4\\\hline
$E^1_{m}$&$a$&$aba$&$abaaaba$&$abaaabababaaaba$\\\hline
$E^2_{m}$&$aa$&$ababa$&$abaaabaaaba$&$abaaabababaaabababaaaba$\\\hline
${}^1\!E^1_{m}$&$a$&$abba$&$abbaaabba$&$abbaaabbabbabbaaabba$\\\hline
${}^1\!E^2_{m}$&$b$&$aa$&$abbabba$&$abbaaabbaaabba$\\\hline
${}^2\!E^1_{m}$&$a$&$aba$&$abacacaba$&$abacacabababacacaba$\\\hline
${}^2\!E^2_{m}$&$aca$&$ababa$&$abacacabacacaba$&$abacacabababacacabababacacaba$\\\hline
\end{tabular}
\end{center}

\medskip

Tab.4 and Tab.5 give the return word sequences of ${}^i\!E^{j}_{m}$ in sequence $\T_i$ for $i,j=1,2$.
Tab.6 gives examples for six envelope words:
$E^1_1$, $E^2_1$, ${}^1\!E^1_1=\tau_1(E^1_1)$, ${}^1\!E^2_2=\tau_1(E^2_1)$,
${}^2\!E^1_2=\tau_2(E^1_1)a$ and ${}^2\!E^2_2=\tau_2(E^2_1)a$.
In these tables, we write ``the return word sequence" as \emph{r.w.s} for short.

For instance, (1) the return word sequence $\{r_p({}^1\!E^{1}_{m})\}_{p\geq1}$ in $\T_1$ is still $\T_1$ over alphabet
$\{r_1({}^1\!E^{1}_{m}),r_2({}^1\!E^{1}_{m})\}=\{\tau_1(r_1(E^{1}_{m})),\tau_1(r_2(E^{1}_{m}))\}$, see Tab.4;
(2) the return word sequence $\{r_p({}^1\!E^1_{1})\}_{p\geq1}=\{r_p(a)\}_{p\geq1}$ in $\T_1$ is still $\T_1$ over alphabet
$\{r_1(a),r_2(a)\}=\{\tau_1(ab),\tau_1(a)\}=\{abb,a\}$, see Tab.6.

\medskip

\small
\begin{center}
\begin{tabular}{|l|c|l|l|l|}
\multicolumn{5}{c}{\normalsize Tab.4: The return word sequences of ${}^1\!E^{j}_{m}$ in sequence $\T_1$ for $j=1,2$}\\\hline
&&\multicolumn{3}{c|}{Alphabet}\\\cline{3-5}
Factor $\w$&r.w.s&$r_1(\w)$&$r_2(\w)$&$r_4(\w)$\\\hline
${}^1\!E^{1}_{m}=\tau_1(E^1_{m})$&$\T_1$
&$\tau_1(r_1(E^{1}_{m}))$&$\tau_1(r_2(E^{1}_{m}))$&$--$\\\hline
${}^1\!E^{2}_{1}=b$&$\T_2$
&$b$&$baaa$&$ba$\\\hline
${}^1\!E^{2}_{m}=\tau_1(E^2_{m-1}),m\geq2$&$\T_2$
&$\tau_1(r_1(E^{2}_{m-1}))$&$\tau_1(r_2(E^{2}_{m-1}))$&$\tau_1(r_4(E^{2}_{m-1}))$\\\hline
\end{tabular}
\end{center}

\begin{center}
\begin{tabular}{|l|c|l|l|l|}
\multicolumn{5}{c}{\normalsize Tab.5: The return word sequences of ${}^2\!E^{j}_{m}$ in sequence $\T_2$ for $j=1,2$}\\\hline
&&\multicolumn{3}{c|}{Alphabet}\\\cline{3-5}
Factor $\w$&r.w.s&$r_1(\w)$&$r_2(\w)$&$r_4(\w)$\\\hline
${}^2\!E^1_{1}=a$&$\T_1$&$ab$&$ac$&$--$\\\hline
${}^2\!E^{1}_{m}=\tau_2(E^1_{m-1})a,m\geq2$&$\T_1$
&$\tau_2(r_1(E^{1}_{m-1}))$&$\tau_2(r_2(E^{1}_{m-1}))$&$--$\\\hline
${}^2\!E^2_{1}=aca$&$\T_2$&$ac$&$acababab$&$acab$\\\hline
${}^2\!E^2_{m}=\tau_2(E^2_{m-1})a,m\geq2$&$\T_2$
&$\tau_2(r_1(E^{2}_{m-1}))$&$\tau_2(r_2(E^{2}_{m-1}))$&$\tau_2(r_4(E^{2}_{m-1}))$\\\hline
\end{tabular}
\end{center}

\normalsize

\begin{center}
\begin{tabular}{|l|c|l|l|l|}
\multicolumn{5}{c}{\normalsize Tab.6: Some examples of the return word sequences}\\\hline
&&\multicolumn{3}{c|}{Alphabet}\\\cline{3-5}
Factor $\w$&$r.w.s$&$r_1(\w)$&$r_2(\w)$&$r_4(\w)$\\\hline
$E^1_1=a\prec \D$&$\T_1$&$ab$&$a$&$--$\\\hline
$E^2_{1}=aa \prec \D$&$\T_2$&$a$&$aababab$&$aab$\\\hline
${}^1\!E^1_{1}=a\prec \T_1$&$\T_1$&$abb$&$a$&$--$\\\hline
${}^1\!E^2_{2}=aa\prec \T_1$&$\T_2$&$a$&$aabbabbabb$&$aabb$\\\hline
${}^2\!E^1_{2}=aba\prec \T_2$&$\T_1$&$abacac$&$ab$&$--$\\\hline
${}^2\!E^2_{2}=ababa\prec \T_2$&$\T_2$&$ab$&$ababacacabacacabacac$&$ababacac$\\\hline
\end{tabular}
\end{center}


\section{Application: combinatorics properties of factors}

Let $\PP $ be a property and $\w$ be a factor.
We call $(\w,p)\in\PP$ if the $p$-th occurrence of $\w$ has property $\PP$.
We call $\w\in\PP$ if there exists $p\geq1$ such that $(\w,p)\in\PP$.
For instance, we consider $\PP$ is ``property of adjacent''. Then $(\w,p)\in\PP$ means there exists an integer $q$ $(>p)$ such that the $p$-th and $q$-th occurrences of $\w$ are adjacent.
And $\w\in\PP$ means there exists two integers $p$ and $q$ such that the $p$-th and $q$-th occurrences of $\w$ are adjacent.
Take $\w=ab\in \D$ for example. The first three occurrences of $\w$ are $\D[1,2]$,
$\D[5,6]$ and $\D[7,8]$. By the definition above, $(\w,1)\not\in \PP $, $(\w,2)\in \PP $ and $\w\in \PP$.
By $\w\in \PP$, we know $\w\w=abab$ is a square in $\D$. But we do not know which integer $p$ such that $(\w,p)\in \PP$.

We are led naturally to study two questions for any sequence $\theta$:

\textbf{Q1.} Determine all factors $\w\in \theta$ such that $\w\in \PP$.

\textbf{Q2.} Determine all $\w\in \theta$ and all integers $p$ such that $(\w,p)\in \PP$.

From our knowledge, most of the previous study on combinatorics on words concern Q1.
Using the return word sequences, we can study Q2. It is more difficult, but useful and interesting.
To study Q2, we introduce a new notion called the \emph{spectrum} of property $\PP$,
which considers both the variables $\w$ and $p$.

\smallskip

Now we study some combinatorial properties: separated ($\PP_1$), adjacent ($\PP_2$) and overlapped ($\PP_3$).
For instance, $(\w,p)\in\PP_1$ means
the $p$-th and $(p+1)$-th occurrences of $\w$ are separated.
By $|r_p(\w)|=|r_p(\E(\w))|$ proved in Section 4, and
the values of $|r_p(E^i_m)|$ ($i=1,2$) proved in Section 2, we have:

Case 1. $\E(\w)=E^1_{m}$.
By Lemma~\ref{L3.5}(1), each factor $\w$ with envelope $E^1_{m}$ satisfies $2^{m-2}+1\leq|\w|\leq2^{m}-1$.
By Theorem \ref{T2.1}, $\{r_p(\w)\}_{p\geq1}$ in $\D$ is $\Theta_1$.
When $\Theta_1[p]=a$,
$|r_p(\w)|=|r_1(\w)|=2^{m}>|\w|$. So the $p$-th and $(p+1)$-th occurrences of $\w$ are separated.
But when $\Theta_1[p]=b$, $|r_p(\w)|=|r_2(\w)|=2^{m-1}$.
In this case, $(\w,p)$ has different property for different $|\w|$.

Case 2. $\E(\w)=E^2_{m}$. By Lemma~\ref{L3.5}(2), each factor $\w$ with envelope $E^2_{m}$ satisfies $2^{m-1}+1\leq|\w|\leq3\times2^{m-1}-1.$
By Theorem \ref{T2.2}, $\{r_p(\w)\}_{p\geq1}$ in $\D$ is $\Theta_2$.
(i) When $\T_2[p]=a$, $|r_p(\w)|=|r_1(\w)|=2^{m-1}<|\w|$;
(ii) when $\T_2[p]=b$, $|r_p(\w)|=|r_2(\w)|=7\times2^{m-1}>|\w|$;
(iii) when $\T_2[p]=c$, $|r_p(\w)|=|r_4(\w)|=3\times2^{m-1}>|\w|$.

Moreover, notice that $\big\{|r_p(\w)r_{p+1}(\w)|,p\geq1\big\}$ is equal to

(1) $\big\{|r_1(\w)r_1(\w)|,|r_1(\w)r_2(\w)|,|r_2(\w)r_2(\w)|\big\}
=\{2^{m+1},3\times 2^m,2^m\}$ in case 1;

(2) $\big\{|r_1(\w)r_2(\w)|,|r_1(\w)r_4(\w)|\big\}
=\{2^{m+2},2^{m+1}\}$ in case 2.

Thus the $p$-th and $q$-th occurrences of $\w$ are always separated for $q-p\geq2$.
So we only discuss the combinatorial properties between the $p$-th and $(p+1)$-th occurrences of a fixed factor.

\smallskip

\textbf{The spectrum of properties $\PP_i$ for $i=1,2,3$.}
\begin{equation*}
\begin{split}
(\w,p)\in\PP_1
\Leftrightarrow~&(\w,p)\in\begin{cases}
\quad\{\E(\w)=E^1_m(m\geq1),\Theta_1[p]=a\}\\
\cup~\{\E(\w)=E^1_m(m\geq1),\Theta_1[p]=b,|\w|<2^{m-1}\}\\
\cup~\{\E(\w)=E^2_m(m\geq1),\Theta_2[p]\neq a\}.
\end{cases}\\
(\w,p)\in\PP_2
\Leftrightarrow~&(\w,p)\in\{\E(\w)=E^1_m(m\geq1),\Theta_1[p]=b,|\w|=2^{m-1}\}\\
(\w,p)\in\PP_3
\Leftrightarrow~&(\w,p)\in\begin{cases}
\quad\{\E(\w)=E^1_m(m\geq1),\Theta_1[p]=b,|\w|>2^{m-1}\}\\
\cup~\{\E(\w)=E^2_m(m\geq1),\Theta_2[p]=a\}.
\end{cases}
\end{split}
\end{equation*}

These are answers of Q2. Using them we have
$$\w\in\PP_2\Leftrightarrow \w\in\{\E(\w)=E^1_m(m\geq1),|\w|=2^{m-1}\},$$
which is an answer of property $\PP_2$ in Q1.
Using this conclusion, we can determine all occurrences of squares and cubes in $\D$.


\end{document}